\documentclass[letterpaper,aps, 12pt,reqno]{amsart}
\usepackage{amsmath,amsthm,amssymb,amsbsy,amsfonts,mathrsfs}
\usepackage{microtype}
\usepackage{a4wide}

\newtheorem{thm}{Theorem}[section]
\newtheorem{prop}[thm]{Proposition}
\newtheorem{lem}[thm]{Lemma}
\newtheorem{cor}[thm]{Corollary}
\theoremstyle{definition}

\theoremstyle{remark}
\newtheorem{remark}[thm]{Remark}
\numberwithin{equation}{section}
\catcode`@=11
\def\section{\@startsection{section}{1}%
  \z@{1.5\linespacing\@plus\linespacing}{.5\linespacing}%
  {\normalfont\bfseries\large\centering}}
\catcode`@=12
\begin{document}

\title[Solitary waves for multicomponent LW-SW system]{
Existence and Positivity properties of Solitary waves for a multicomponent
long wave--short wave Interaction system}
\author[Santosh Bhattarai]{Santosh Bhattarai}
\address{Trocaire College,
360 Choate Ave, Buffalo, NY 14220 USA}
\email{sntbhattarai@gmail.com, bhattarais@trocaire.edu}


\begin{abstract}
We study the existence of solitary-wave solutions and some of their properties for
a general multicomponent long-wave--short-wave interaction system.
The system considered here describes the nonlinear interaction of multiple short waves with a long-wave, and
is of interest in plasma physics, nonlinear optics, and fluid dynamics.
\end{abstract}
 \maketitle
\section{Introduction}

\noindent The long wave-short wave interaction (LSI) is an important problem in a variety of physical systems.
The LSI model has been successfully
applied to many different contexts of modern physics and fluid dynamics, such as
studying the solitons resulting from the interactions between long ion-sound waves (ion-acoustic waves)
and short Langmuir waves (plasma waves, plasmons) in a magnetized
plasma \cite{[Kar],[YO]}, or Alfv\'{e}n-magneto-acoustic waves interactions in a cold plasma \cite{[Sul]}.
Kawahara et al. \cite{[KSK]} have investigated the nonlinear interaction
between short and long capillary-gravity waves on a liquid layer of uniform depth.
For a general theory for deriving nonlinear PDEs which permit both
long and short wave solutions and interact each other nonlinearly, the reader may consult \cite{[Ben]}.

\smallskip

\noindent In recent years there has been renewed interest in the study of nonlinear waves in multi-component LSI system.
The multi-component LSI systems arise in water waves theory \cite{[Cra]}, optics \cite{[Ohta1]},
ferromagnetism theory \cite{[My]}, acoustics \cite{[Saz1]}, in a bulk elastic medium \cite{[Erb]}, to name a few.
In this paper we consider a general multi-component LSI system
describing the interaction of multiple NLS-type short waves with a KdV-type long wave in one dimension.
The nonlinear interaction between $N$ complex short-wave field envelopes, $u_j, j=1,2,...,N,$ and the real long-wave field, $v,$ can be modeled by the $(N+1)$-component long-wave--short-wave system
\begin{equation}\label{lwsw}
\left\{
\begin{aligned}
  & i\partial_t u_1 + \partial_{x}^{2} u_1 = - \alpha_1 u_1v, \\
  & \ldots \ \ \ \ \ \  \ldots \  \ \ \ \  \ \ \ \ \ldots  \\
  & i\partial_t u_N + \partial_{x}^{2} u_N = - \alpha_N u_Nv,\\
   & \partial_t v+\partial_x\left(\gamma \partial_x^{2}v+Q(v)\right) = - \partial_{x} \left(\beta_1|u_1|^2+\ldots +\beta_N|u_N|^2 \right),\\
\end{aligned}
\right.
\end{equation}
where $x$ and $t$ are spatial and temporal variables, respectively, $Q=Q(v)$ is a nonlinear polynomial,
and the parameters $\alpha_j, \beta_j, \gamma$ are real constants.
The motivation for studying systems of the form \eqref{lwsw} also come from a pioneer work of Kanna et al.~\cite{[Kanna]}, who
set $\gamma=0, \alpha_j=\alpha=\beta_j,$ and $Q\equiv0$ in \eqref{lwsw} and proved that
\eqref{lwsw} can be derived from a system of multi-component coupled nonlinear Schr\"{o}dinger type equations.
They have also shown that the system is integrable via Painlev\'{e} test.
In \cite{[Ma1981]}, the system \eqref{lwsw} with $\gamma=0, \alpha_j=\alpha=\beta_j, Q\equiv0,$ and $N=2$ has been
shown to be integrable by the inverse scattering transform method and the soliton solutions have been obtained. In the same case,
the rogue waves of \eqref{lwsw} have been reported in \cite{[Ch1]}. We also mention the
paper \cite{[Ch2]} where general bright-dark multi-soliton solution has been constructed
for a general multicomponent LSI system.

\medskip

\noindent The mathematical study of
systems of the form \eqref{lwsw} with $N=1$ and $Q(v)=v^2,$ namely well-posedness
results (unique existence, persistence, and continuous dependence on initial data)
on the associated Cauchy problem or existence and qualitative properties of solutions, has been
studied extensively over the years by many authors using both numerical and theoretical techniques
(see for example, \cite{[AP],[Bai1],[Cha1],[Colo],WU} and references therein).
Despite some progress has been made so far using numerical and algorithms methods,
many difficult questions remain open and
little is known about theoretical results concerning existence
and properties of solutions for $(N+1)$-systems \eqref{lwsw} for $N\geq 2.$
This project aims to cast a light on $(N+1)$-component long-wave--short-wave interaction
system. Included in the study are existence results and several
properties of travelling solitary waves for \eqref{lwsw} in
the case when $\beta_j=\alpha_j/2$ and $Q(v)=\tau v +\beta v^2$ with $\tau \in \mathbb{R} $ and $\beta\geq 0.$

\medskip

\noindent By the travelling-wave solutions of \eqref{lwsw} we mean the solutions of the form
\begin{equation}\label{SO}
\mathbb{T}=\left\{ \left(e^{i\omega t}\phi_1(x-ct),\ldots , e^{i\omega t} \phi_N(x-ct), \Psi(x-ct)\right):c, \omega\in \mathbb{R} \right\}.
\end{equation}
Usually a nontrivial (i.e., not identically zero) travelling-wave
solution which vanishes at $\pm \infty$ (say, $\phi_j$ and $\Psi$ are in $H^1(\mathbb{R}),$ the usual Sobolev space)
is referred to as a solitary wave.
In the case when $c=0$ (zero travelling velocity), these solutions \eqref{SO} are time independent which usually are
referred to as standing-wave solutions or stationary-state solutions. Let $(u_1,\ldots, u_N, v)$ be
a solution of the form \eqref{SO}. Put $\phi_j(x)=e^{icx/2}\Phi_j(x), j=1,\ldots, N,$ and
substitute $u_j$ and $v$ into \eqref{lwsw}, integrate the second equation once,
and evaluate the constant of
integration by using the fact that $\Phi_j$ and $\Psi$ are $H^1$ functions.
Then, one sees that $(\Phi_1,\ldots,\Phi_N,\Psi)$
must satisfy
the following system of ordinary
differential equations
\begin{equation}\label{ODE}
\left\{
\begin{aligned}
& \Phi_1^{\prime \prime}-\sigma \Phi_1=-\alpha_1 \Phi_1 \Psi , \\
& \ldots \ \ \ \ \ \ldots \  \ \ \ \ \ldots \\
& \Phi_N^{\prime \prime}-\sigma \Phi_N=-\alpha_N \Phi_N \Psi, \\
& \gamma\Psi^{\prime \prime}-c_\tau\Psi =- \frac{1}{2}\beta \Psi^2- \frac{1}{2}\left(\alpha_1\Phi_j^{2}+\ldots + \alpha_N\Phi_N^{2}\right),
\end{aligned}
\right.
\end{equation}
where $c_\tau=c-\tau,\ \sigma=\omega-c^2/4,$ and primes
denote derivatives with respect to the variable $\xi.$
This equation will be considered the defining equation for solitary waves.

\medskip

\noindent Our analysis begins with a study of the existence of solutions of the system \eqref{ODE}.
We prove the existence of nontrivial
solutions of \eqref{ODE} with each component in $H^\infty$ and exponential decay at infinity under the
following assumptions
\begin{equation}\label{assumptions}
c>0, \ \sigma>0,\  -\infty < \tau \leq c,\ \beta\geq 0,\ \gamma\geq 0,\ \alpha_j>0
\end{equation}
(see Theorem~\ref{existence} below
for the precise statement of the result). The existence result is proved by studying a minimization problem whose
minimizers, up to a constant, corresponds to solitary waves for \eqref{lwsw}.
More precisely, let $\mathsf{K}$ be the functional defined for $(f_1,...,f_N,g)\in (H^1)^{N+1}$ by
\begin{equation}
\mathsf{K}(f_1, \ldots, f_N,g)=\int_{-\infty}^{\infty}(f_1, \ldots, f_N,g)D(L_{ii})(f_1, \ldots, f_N,g)^T\ dx,
\end{equation}
where $D(L_{ii})$ is the diagonal matrix with diagonal entries $L_{ii}=\sigma-\partial_{xx}$ for $i=1,\ldots, N,$ and $L_{ii}=c_\tau-\gamma\partial_{xx}$ for $i=N+1,$ and
introduce the notation
\begin{equation*}
F(f_1,\ldots,f_N,g)=\frac{1}{3} \beta g^3+\left( \alpha_1 f_{1}^{2}+\alpha_2 f_{2}^{2}+\ldots + \alpha_N f_{N}^{2}\right)g.
\end{equation*}
For $\lambda >0,$ we shall show that the variational problem (P1) of minimizing the
functional $\mathsf{K}(f_1,...,f_N,g)$ subject to the constraint
\begin{equation*}
\int_{-\infty}^{\infty}F(f_1,\ldots,f_N,g)(x)\ dx =\lambda
\end{equation*}
always has a non-empty solution set provided that \eqref{assumptions} holds.
The key idea in establishing the existence of minimizers
here is to apply Lions' concentration compactness
lemma (Lemma~\ref{concentration}) to a minimizing sequence of the problem (P1)
and extract a subsequence which is tight. The method of concentration compactness then implies that
this subsequence, when its terms are suitably translated, converges
strongly in $(H^1)^{N+1}$ to a limit which achieves the minimum of the problem (P1).
Let $(\tilde{\phi}_1,...,\tilde{\phi}_N,\tilde{\psi})$ be this limit. Then, by the Lagrange multiplier principle,
this limit function $(\tilde{\phi}_1,...,\tilde{\phi}_N,\tilde{\psi}),$ after
multiplying by a constant, corresponds to a solution of \eqref{ODE}, at least in the distributional sense.
Such solutions are called weak ground state solutions.
But since the right sides of all $N+1$ equations in the system \eqref{ODE} are continuous functions,
a standard bootstrapping argument shows that weak ground state solutions are
indeed classical solutions (see Proposition~\ref{smoothdecay} below).
In Section 2, we provide the details of the method.

\medskip

\noindent In Section 3, we combine the variational
formulation of solutions of \eqref{ODE} with the theory of
symmetric decreasing rearrangements to prove the existence of
solutions $(\Phi_j,\Psi)$ such that $\Phi_j$ and $\Psi$ are even and decreasing
positive functions in $(0,\infty)$ (see Theorems~\ref{signthm} and \ref{evenSOL} below).
Similar techniques have been used previously by Albert et al \cite{[AB2]} to study
solitary-wave solutions of some model equations for waves in stratified fluids, and by
Angulo and Montenegro \cite{[AP]} to prove the existence and evenness of
solitary waves for an interaction equation in two-layer fluid.
The paper closes by showing in Theorem~\ref{positiveFT} the existence of a
solitary wave for \eqref{lwsw} with positive Fourier transforms.
We state these results in the context of the existence theory for solitary waves
introduced by Weinstein in \cite{[Wein]}
and will be proved by adapting an argument developed by Albert in \cite{[Albp]}.

\medskip

\noindent \textbf{Notation.}
For $(x, a)\in \mathbb{R}\times (0,\infty),$ we denote by $B(x,a)$ the ball centered at $x$ and of radius $a.$
In particular, we denote $B_a=B(0,a).$
The Fourier transform $\widehat{f}$ of a tempered distribution $f(x)$ on $\mathbb{R}$ is defined as
$\widehat{f}(\xi)=(2\pi)^{-1/2}\int_{-\infty}^{\infty}f(x)e^{i\xi x}\ dx.$ If $1\leq r<\infty,$ we shall
denote by $L^r=L^r(\mathbb{R})$ the usual
Banach  space  of Lebesgue measurable
functions $f$ on $\mathbb{R}$ for which the norm $|f|_{L^r}$ is finite, where
\begin{equation*}
|f|_{L^r}=\left( \int_{-\infty }^{\infty }\left\vert
f\right\vert ^{r}dx\right) ^{1/r}\textrm{ \ for }\ 1\leq r<\infty.
\end{equation*}
The space $L^\infty$ consists of the measurable, essentially bounded functions on $\mathbb{R}$ with the norm
$|f| _{L^\infty }=\textrm{ess} \sup_{x\in \mathbb{R}}|f(x)|.$
The (Lebesgue) convolution of two functions $f$ and $g,$ denoted by $f\star g,$ is the integral
\begin{equation*}
f\star g(x)=f(x)\star g(x)=\int_{-\infty}^{\infty}f(x-\xi)g(\xi)\ d\xi.
\end{equation*}
For $s\geq 0,$ we denote by $H^s=H^{s}(\mathbb{R})$ the Sobolev space of all tempered
distributions $f$ on $\mathbb{R}$ whose Fourier transforms $\widehat{f}$ are measurable functions on $\mathbb{R}$ satisfying
\begin{equation*}
\| f\| _{s}^{2}=\int_{-\infty }^{\infty }\left(
1+|\xi|^{2}+...+|\xi|^{2s}\right)|\widehat{f}(\xi
)| ^{2}d\xi <\infty .
\end{equation*}
In particular, we use $\|
f\|$ to denote the $L^{2}$ or $H^{0}$ norm of a
function $f.$
We define the space $\mathcal{Y}$ to be the cartesian product $H^1 \times . . . \times H^1$ ($N+1$-times) provided
with the product norm $\|\cdot\|_{\mathcal{Y}}.$
For notational convenience, we denote
\begin{equation*}
\begin{aligned}
& (\mathbf{u}_j,v)=(u_1,...,u_N,v),\\
& (\mathbf{u}_{j,n},v_n)=(u_{1,n},...,u_{N,n},v_n),\ \textrm{and}\
 (\Phi_j,\Psi)=(\Phi_1,...,\Phi_N,\Psi).
\end{aligned}
\end{equation*}
In place of the compound subscripts, for example, when we
take a subsequence of a sequence, we will follow the convention of
using the same symbol to denote the subsequence. The letter $C$ will be used to
denote various positive constants which may assume different values from line to
line but are not essential
to the analysis of the problem. The letter $C(...)$ will denote the constant whose value depends essentially only
on the quantities indicated in the parentheses.

\section{Existence of Solitary Waves}
\noindent The main result of this section is the existence of global minimizers for the variational problem (P1):

\begin{thm}\label{existence}
 Suppose that the assumptions \eqref{assumptions} hold for the
 constants $c,\ \sigma,\ \beta,\ \tau, \gamma,$ and $\alpha_j.$ For $\lambda>0,$ define
 \begin{equation*}
 \mathcal{A}=\left\{(f_1,...,f_N,g)\in \mathcal{Y}:\int_{-\infty}^{\infty} F(f_1,...,f_N,g)(x)\ dx =\lambda \right\}.
 \end{equation*}
Then there exists a minimizing function for the problem (P1) in $\mathcal{A}.$ Consequently, the system \eqref{ODE} has
 a solution $(\Phi_1,...,\Phi_N,\Psi)$ such
 that $\Phi_1,...,\Phi_N, \Psi$ are in $H^{\infty}(\mathbb{R})$ and decay exponentially at infinity.
\end{thm}
\noindent In particular, Theorem~\ref{existence} guarantees that
the minimizing set $S(\lambda),$ namely
\begin{equation*}
S(\lambda)=\left\{(\Phi_j,\Psi)\in \mathcal{A}:\mathsf{K}(\Phi_j,\Psi)=\inf\mathsf{K}(\mathbf{f}_j,g),\ (\mathbf{f}_j,g)\in \mathcal{A}\right\},
\end{equation*}
 is non-empty. As will be seen below, this translates into an
existence result for solitary-wave solutions \eqref{SO} of \eqref{lwsw}.

\smallskip

\noindent We begin by proving some properties of the variational problem.
The first lemma shows that $\mathsf{K}$ has a finite and positive infimum on $\mathcal{A}.$

\begin{lem}\label{monotono}
For each $\lambda>0,$ one has
\begin{equation}\label{Idef}
I_\lambda=\inf\left\{\mathsf{K}(f_1,...,f_N,g):(f_1,...,f_N,g)\in \mathcal{A}\right\}>0.
\end{equation}
Moreover, if $\lambda_2>\lambda_1>0,$ then $I_{\lambda_2}\geq I_{\lambda_1}.$
\end{lem}
\noindent{\bf Proof.}
Denote $\Delta=(f_1,...,f_N,g).$ From the Cauchy-Schwartz
inequality and the Sobolev embedding theorem we have
\begin{equation}\label{Ipositive}
\begin{aligned}
\lambda & =\int_{-\infty}^{\infty}F(\Delta)(x)\ dx
& \leq C\left( \|g\|_1 \|g\|^2+ \sum_{j=1}^{N}\|f_j\|_1 \|f_j\| \|g\| \right)\leq C\|\Delta\|_{\mathcal{Y}}^{3},
\end{aligned}
\end{equation}
where the constant $C$ is independent of $f_j, \ 1\leq j \leq N,$ and $g.$ Now, using \eqref{Ipositive} it follows that
\begin{equation*}
\begin{aligned}
\mathsf{K}(\Delta)& \geq \min\{1,\sigma\}\sum_{j=1}^{N} \|f_j\|_{1}^{2}+\min\{\gamma,c_\tau \}\|g \|_{1}^{2}\\
& \geq \min\{\min\{1,\sigma\}, \min\{\gamma,c_\tau \}\} \|\Delta\|_{\mathcal{Y}}^{2}\geq C\lambda^{2/3}>0,
\end{aligned}
\end{equation*}
and therefore $I_\lambda>0.$ To prove $I_{\lambda_2}\geq I_{\lambda_1},$
let $\epsilon>0$ be arbitrary. There exists a function $\Theta=(\phi_1,...,\phi_N,\psi)$ in $\mathcal{Y}$
such that $\int_{-\infty}^{\infty}F(\Theta)dx=\lambda_2$ and
$\mathsf{K}(\Theta)<I_{\lambda_2}+\epsilon.$ For $a\in \mathbb{R},$ denote
\begin{equation*}
Q(a\Theta)=\int_{-\infty}^{\infty}F(\Theta(x))\ dx.\
\end{equation*}
Then $Q(a\Theta)$ is a continuous function of $a\in \mathbb{R}$ and hence, using the intermediate
value theorem of elementary analysis, we can find $\xi\in (0,1)$ such that $Q(\xi\Theta)=\lambda_1.$ Hence
\begin{equation*}
I_{\lambda_1}\leq \mathsf{K}(\xi\Theta)=\xi^2\mathsf{K}(\Theta)<\mathsf{K}(\Theta)<I_{\lambda_2}+\epsilon.
\end{equation*}
Since $\epsilon>0$ is arbitrary, it follows that $I_{\lambda_1}\leq I_{\lambda_2},$ proving the lemma. \hfill $\Box$

\medskip

\noindent By a minimizing sequence for $I_{\lambda }$ in what follows, we mean to be any sequence $
\{(\mathbf{f}_{j, n},g_{n})\}$ of functions in $\mathcal{A}$ satisfying the conditions
\begin{equation}\label{minimizingseq}
\lim_{n\to \infty }\mathsf{K}(\mathbf{f}_{j, n},g_{n})=I_{\lambda }\
 \textrm{and} \  \int_{-\infty}^{\infty} F(\mathbf{f}_{j, n},g_{n})(x)dx=\lambda ,\ \forall n.
\end{equation}

\begin{lem}\label{strictmlambda}
For all $\lambda>0$ and $m>1,$ one has $I_{m \lambda}< mI_{\lambda}.$
\end{lem}
\noindent{\bf Proof.}
Let $\{(\mathbf{f}_{j,n},g_n)\}$ be any sequence of functions in $\mathcal{A}$
satisfying \eqref{minimizingseq}. Denote $\Delta_n=(\mathbf{f}_{j,n},g_n).$ Choose $\theta_n>0$ such that
\begin{equation}\label{strictQ}
\int_{-\infty}^{\infty}F(\theta_n\Delta_{n})(x)\ dx= m\lambda.
\end{equation}
Since $\int_{-\infty}^{\infty}F(\Delta_n) dx =\lambda,$ it follows from \eqref{strictQ} that $\theta_{n}^3=m>1.$ Thus
\begin{equation*}
I_{m\lambda}\leq \mathsf{K}(\theta_n \Delta_{n})=\frac{m}{\theta_n}\mathsf{K}(\Delta_{n}).
\end{equation*}
Since $m>1$ and there exists $\epsilon>0$ such that $\theta_n>1+\epsilon$ for sufficiently large $n,$ the
lemma follows by letting $n\to \infty$ in the last inequality. \hfill $\Box$

\medskip

\noindent As an immediate corollary of Lemma~\ref{strictmlambda}, we obtain the following strict subadditivity property of $I_{\lambda}:$

\begin{cor}\label{strictsubadd}
Let $I_\lambda$ be as defined in \eqref{Idef}.
Then, for all $\lambda_1, \lambda_2>0,$
\begin{equation*}
I_{(\lambda_1+\lambda_2)}< I_{\lambda_1}+I_{\lambda_1}.
\end{equation*}
\end{cor}
\noindent{\bf Proof.}
Without loss of generality, we may assume that $\lambda_1\geq \lambda_2.$ If $\lambda_1> \lambda_2,$ then from what was shown in
Lemma~\ref{strictmlambda}, it follows that
\begin{equation*}
\begin{aligned}
I_{(\lambda_1+\lambda_2)}& = I_{\lambda_1\left(1+\lambda_2 \lambda_{1}^{-1}\right)}
<\left(1+\lambda_2\lambda_{1}^{-1}\right)I_{\lambda_1} \\
& \leq I_{\lambda_1}+\lambda_2\lambda_{1}^{-1} \left(\lambda_1 \lambda_{2}^{-1} I_{\lambda_2}\right) = I_{\lambda_1}+I_{\lambda_2};
\end{aligned}
\end{equation*}
whereas in the case $\lambda_1= \lambda_2,$ we have
\begin{equation*}
I_{(\lambda_1+\lambda_2)}=I_{2\lambda_1}<2I_{\lambda_1}=I_{\lambda_1}+I_{\lambda_2},
\end{equation*}
so the corollary has been proved. \hfill $\Box$

\medskip

\noindent We now proceed to prove the existence result.
The idea is to show that any
minimizing sequence
$\{(\mathbf{f}_{j, n}, g_{n})\}_{n\in \mathbb{N}}$
for $I_\lambda$ in $\mathcal{A}$ which, up to subsequences and when its
terms are suitably translated, has the following properties:
\begin{equation}\label{solved1}
(\Phi_j,\Psi)=\lim_{n\to \infty}\Delta_n \in \mathcal{A}\ \ \textrm{and} \ \mathsf{K}(\Delta_n)\leq \liminf_{n\to \infty}\mathsf{K}(\Delta_n),
\end{equation}
where $\Delta_n=(\mathbf{f}_{j, n}, g_{n}).$ Once we establish \eqref{solved1}, the minimization
problem (P1) is then solved since then it follows that
\begin{equation*}
I_\lambda\leq \mathsf{K}(\Phi_j,\Psi)\leq \liminf_{n\to \infty}\mathsf{K}(\Delta_n)=I_\lambda,
\end{equation*}
where the first inequality holds because the limit pair $(\Phi_j,\Psi)$ belongs to $\mathcal{A}$ and the second
inequality holds because $\{\Delta_n\}_{n\in \mathbb{N}}$ is a minimizing sequence for $I_\lambda.$
The key tool here is the concentration compactness
principle developed by P.~L.~Lions \cite{[L1]}, which has been used
by many authors (see, for example, \cite{[A],[AB2],[AP],[San3],[SB11],[CBona], [Zeng]} and references therein).
The method is based on the following lemma:

\begin{lem}[Lions \cite{[L1]}]\label{concentration}
Let $\left\{Q _{n}\right\} _{n\geq 1}$ be a sequence of nonnegative
functions in $L^{1}(\mathbb{R})$ satisfying $\int_{-\infty }^{\infty }Q
_{n}(x)\ dx=\alpha $ for all $n$ and some fixed $\alpha >0.$ Then there exists a
subsequence $\left\{ Q _{n_{k}}\right\} _{k\geq 1}$ satisfying exactly one of the
following three possibilities:
\begin{itemize}
\item[(1)] (Tightness up to translation) There are $y_{k}\in \mathbb{R}$ for $k=1,2,.\ .\ .,$
such that $Q _{n_{k}}(.+y_{k})$ is tight, i.e., for any $\varepsilon >0,$
there is $R>0$ large enough such that%
\begin{equation*}
\int_{B(y_k, R)}Q _{n_{k}}(x)\ dx\geq \alpha
-\epsilon \ \ \textrm{for all} \ k.
\end{equation*}

\item[(2)] (Vanishing) For any $R>0,$%
\begin{equation*}
\lim_{k\longrightarrow \infty }\sup_{y\in \mathbb{R}}\int_{B(y,R)}Q _{n_{k}}(x)\ dx=0.
\end{equation*}

\item[(3)] (Dichotomy) There exists $\bar{\alpha} \in (0,\alpha )$ such that for any
$\varepsilon >0,$ there exists $k_{0}\geq 1$ and $Q _{1, k},Q
_{2, k} \in L^{1}_{+}(\mathbb{R})$
such that for $k\geq k_{0},$%
\begin{equation*}
\left\{
\begin{aligned}
&\left\vert Q _{n_{k}}-(Q _{1, k}+Q _{2, k})\right\vert _{1}\leq
\epsilon ,\ \ \left\vert \int_{\mathbb{R}}Q _{1, k}\
dx-\bar{\alpha} \right\vert \leq \epsilon , \\
&\left\vert \int_{\mathbb{R} }Q _{2, k}\ dx-(\alpha -\bar{\alpha} )\right\vert \leq \varepsilon ,\\
& \ dist(\textrm{supp}(Q
_{1, k}),\textrm{supp}(Q _{2, k}))\rightarrow \infty \text{ \ as }%
k\rightarrow \infty .%
\end{aligned}
\right.
\end{equation*}
\end{itemize}
\end{lem}

\begin{remark}
In Lemma~\ref{concentration} above, the condition $\int_{\mathbb{R}}Q _{n}(x)\
dx=\alpha $ can be replaced by $\int_{\mathbb{R} }Q _{n}(x)\
dx=\alpha _{n}$ where $\alpha _{n}\to \alpha >0$ as $%
n\to \infty .$ Indeed, it is enough to replace $Q _{n}$ by $%
Q _{n}/\alpha _{n}$ and apply the lemma.
\end{remark}

\noindent We now consider a minimizing sequence $\{(\mathbf{f}_{j, n},g_{n})\}_{n\in \mathbb{N}}$ for $I_\lambda$
and apply the Lemma~\ref{concentration}. Denote $\Delta_n=(\mathbf{f}_{j, n},g_{n})$ and let
\begin{equation}\label{Qdef}
Q _{n}=\left(g_{n}^{\prime}\right)^2+g_{n}^{2}+ \sum_{j=1}^{N}\left((f_{j, n}^{\prime
})^{2}+ (f_{j, n})^{2}\right).
\end{equation}
For each $n,$ define $\mu_n=\int_{-\infty}^{\infty}Q _{n}(x)\ dx.$
As $\{\Delta_n\}_{n\in \mathbb{N}}$ is a minimizing sequence,
the sequence $\{\mu_n \}_{n\in \mathbb{N}}$ of real numbers is uniformly bounded for sufficiently large $n.$
Without loss of generality, suppose that $\int_{-\infty}^{\infty}Q _{n}(x)\ dx \to \mu$ whenever $n\to \infty.$
By Lemma~\ref{concentration} above, the sequence $\{Q_n\}_{n\in \mathbb{N}}$ has a subsequence
which satisfies one of the three possibilities: Tightness up to translation, Vanishing, or Dichotomy.
Our task is to show that tightness up to translation is the only possibility. To this end, suppose
there is a subsequence $\{Q_{n_k}\}_{n\in \mathbb{N}}$ of $\{Q_{n} \}_{n\in \mathbb{N}}$ which
satisfies either vanishing or dichotomy. We divide the proof into a sequence of lemmas.
The first lemma rules out the vanishing condition:

\begin{lem}
Vanishing does not occur.
\end{lem}
\noindent{\bf Proof.}
We prove this lemma in several steps.

\smallskip

\noindent \textit{Step 1.} Suppose $g\in C^{\infty}(\mathbb{R}),$ and for $m\in \mathbb{Z},$
define $I_m=[m-1/2, m+1/2].$ Then for all $m\in \mathbb{Z},$ one has
\begin{equation}\label{step1inq}
\sup_{x\in I_m}|g(x)|\leq \int_{I_m}|g(y)|\ dx +\int_{I_m} |g^\prime(y)|\ dy.
\end{equation}
To see this, for all $z\in I_m$ and $y\in I_m,$ it is obvious that
\begin{equation*}
g(z)=g(y)+\int_{y}^{z}g^\prime(x)\ dx.
\end{equation*}
In consequence, one has for all $m\in \mathbb{Z},$
\begin{equation*}
|g(z)|\leq |g(y)|+\int_{I_m}|g^\prime(x)|\ dx.
\end{equation*}
Integrating both sides with respect to $y$ over $I_m,$ one obtains that
\begin{equation*}
|g(z)|\leq \int_{I_m} |g(y)|\ dy +\int_{I_m}|g^\prime(y)|\ dy,
\end{equation*}
from which \eqref{step1inq} follows.

\smallskip

\noindent \textit{Step 2.} Suppose $\Delta=(f_1,\ldots,f_n,g)\in \mathcal{Y}$ and $Q=Q_n$ be as
defined in \eqref{Qdef} with $\Delta_n$ replaced by the constant sequence $\Delta=\Delta_n.$ Then there exists $C>0$ such that
\begin{equation}\label{step3inq1}
\int_{-\infty}^{\infty}|g|^3\ dx \leq C \left( \sup_{y\in \mathbb{R}}\int_{B(y,1/2)}Q(x)\ dx\right)^{1/2} \|\Delta\|_{\mathcal{Y}}^{2}
\end{equation}
and for all $j=1,\ldots,N,$
\begin{equation}\label{step3inq2}
\begin{aligned}
\ \ \ \ \ \ \ \ \int_{-\infty}^{\infty}|g||f_j|^2\ dx \leq C \left( \sup_{y\in \mathbb{R}}\int_{B(y,1/2)}Q(x)\ dx\right)^{1/2} \|\Delta\|_{\mathcal{Y}}^{2}.
\end{aligned}
\end{equation}
To prove \eqref{step3inq1}, assume first that $g\in C_0^{\infty}(\mathbb{R})$. By Step 1, replacing $g$ by $g^2,$ we obtain that
\begin{equation*}
\begin{aligned}
& \left( \sup_{y\in I_m}| g(y)| \right) ^{2} \leq \|g^2(y)\|_{L^1(I_m)}+2\int_{I_m}|g(y)||g^\prime(y)|\ dy\\
&\ \ \ \ \ \ \ \ \leq 2\|g(y) \|_{L^2(I_m)}^2 + \|g^\prime(y)\|_{L^2(I_m)}^2
\leq C\sup_{y\in \mathbb{R}}\int_{B(y,1/2)}Q(x)\ dx.
\end{aligned}
\end{equation*}
Since $\int_{I_m}|g|^3\ dx\leq \|g \|_{L^\infty(I_m)} \|g\|_{L^2(I_m)}^2,$ using the above
estimate and taking the sum over all $m\in \mathbb{Z},$ it follows that
\begin{equation*}
\begin{aligned}
\int_{-\infty}^{\infty}|g|^3\ dx
\leq C \left( \sup_{y\in \mathbb{R}}\int_{B(y,1/2)}Q(x)\ dx\right)^{1/2}\|\Delta\|_{\mathcal{Y}}^2.
\end{aligned}
\end{equation*}
This proves \eqref{step3inq1} for $g\in C_0^{\infty}(\mathbb{R}).$ The result for $g\in H^1$ follows by approximating $g$ with
a sequence $\{g_n\}\subset C_0^{\infty}(\mathbb{R})$ such that $g_n\to g$ in $H^1$ norm.
The proof for \eqref{step3inq2} uses the same argument.

\smallskip

\noindent \textit{Step 3.} Suppose now that the vanishing case occurs, which is to say
\begin{equation*}
\lim_{n\to \infty }\sup_{y\in \mathbb{R}}\int_{B(y,1/2)}Q_{n_k}(x)\ dx=0.
\end{equation*}
But then using the estimates obtained in Step 2, it follows that
\begin{equation*}
\begin{aligned}
0<\lambda= &\left\vert \int_{-\infty}^{\infty}F(\Delta_{n_k}) \ dx\right\vert
\leq C \|\Delta_{n_k}\|_{\mathcal{Y}}^{2}\left(\sup_{y\in \mathbb{R}}\int_{B(y,1/2)}Q_{n_k}(x)\ dx \right)^{1/2}\to 0,
\end{aligned}
\end{equation*}
as $n\to \infty,$ which is a contradiction. \hfill $\Box$

\medskip

\noindent The next step in the proof of Theorem~\ref{existence} is to rule out the dichotomy case.
This is dealt with in the next three
lemmas, which represent a simplification and generalization of arguments appeared in \cite{[A],[L1]}.
\begin{lem}\label{Maindicho}
Suppose there is a subsequence $\{Q_{n_k}\}_{k\in \mathbb{N}}$ of $\{Q_{n}\}_{n\in \mathbb{N}}$ such that
the dichotomy alternative of Lemma~\ref{concentration} holds and denote
$\Delta_{n_k}=(\mathbf{f}_{j, n_k}, g_{n_k}).$
Then there exists a real number $\bar{\lambda}=\bar{\lambda}(\epsilon),$ a natural number $n_0,$ and
two sequences of functions
$\{\Delta_{k}^{(1)}\}$ and $\{\Delta_{k}^{(2)}\}$ in $Y$
satisfying
$\Delta_{n_k}=\Delta_{k}^{(1)}+\Delta_{k}^{(2)}$
for all $k$ and for all $k\geq k_0,$
\begin{equation*}
\begin{aligned}
& (i)\ \int_{-\infty}^{\infty} F(\Delta_{k}^{(1)})\ dx -\bar{\lambda}=O(\epsilon),\\
& (ii)\ \int_{-\infty}^{\infty}F(\Delta_{k}^{(2)})\ dx-(\lambda-\bar{\lambda})=O(\epsilon),\\
& (iii)\ \mathsf{K}(\Delta_{n_k})=\mathsf{K}(\Delta_{k}^{(1)})
+\mathsf{K}(\Delta_{k}^{(2)})+O(\epsilon),
\end{aligned}
\end{equation*}
where the constants implied in the notation $O(\epsilon)$ can be chosen independently of $n$ as well as $\epsilon.$
Furthermore, one has
\begin{equation}\label{Kestimate}
\mathsf{K}(\Delta_{k}^{(1)})\geq \bar{\mu} +O(\epsilon) \ \textrm{and}\ \ \mathsf{K}(\Delta_{k}^{(2)})\geq \mu - \bar{\mu} +O(\epsilon),
\end{equation}
where the real number $\bar{\mu}$ is as defined in Lemma~\ref{concentration}.
\end{lem}

\begin{remark}
The lemma says that the subsequence $\{\Delta_{n_k}\}_{k\in \mathbb{N}}$
can be split into two summands which
carry fixed proportions of the constraint and which
are supported far apart spatially that the sum of the values of the functional $\mathsf{K}$ at each summand
does not exceed $\mathsf{K}(\Delta_{n_k}).$
\end{remark}

\medskip

\noindent{\bf Proof.}
If dichotomy case of Lemma~\ref{concentration} occurs, then there exists $\bar{\mu}\in (0,\mu)$ such that for any $\epsilon>0$ there
corresponds $k_0\geq 1$ and $L^1$ functions
$Q_{1,k}, Q_{2,k}\geq 0$ such that for all $k\geq k_0,$
\begin{equation}\label{dichoeq}
\begin{aligned}
&\left\vert Q _{n_{k}}-(Q_{1, k}+Q_{2, k})\right\vert _{1}\leq
\varepsilon , \\
& \left\vert \int_{-\infty}^{\infty}Q_{1, k} \ dx-\bar{\mu} \right\vert \leq \varepsilon ,\ \ \textrm{and} \ \
\left\vert \int_{-\infty}^{\infty}Q_{2, k}\ dx-(\mu -\bar{\mu} )\right\vert \leq \varepsilon.
\end{aligned}
\end{equation}
Moreover, without loss of generality, we may assume that the
supports of the functions $Q_{1, k}$ and $Q_{2, k}$ are separated as follows:
\begin{equation}\label{supporteq}
\begin{aligned}
&\textrm{supp}\ Q_{1, k}\subset (y_k-R_0, y_k+R_0), \\
&\textrm{supp}\ Q_{2, k}\subset(-\infty, y_k-2R_k)\cup (y_k-2R_k, \infty),
\end{aligned}
\end{equation}
for some fixed $R_0>0,$ a sequence of real numbers $\{y_k\}_{n\in \mathbb{N}},$ and $R_k\to \infty.$
To split $\Delta_{n_k}$ into two summands $\Delta_{k}^{(1)}$ and $\Delta_{k}^{(2)}, \ k=1,2,\ ...,$ let
$\zeta$ and $\rho \in C_{0}^{\infty}(\mathbb{R})$ with $0 \leq \zeta, \rho \leq 1$ be such that
$\zeta\equiv 1$ on  $B_1,$\ $\textrm{supp}\ \zeta \subset B_2 ; \
\rho \equiv 1 \ \textrm{on}\ \mathbb{R}\setminus B_2,\ \textrm{supp}\ \rho \subset \mathbb{R}\setminus B_1.$
Denote the functions
\begin{equation*}
\zeta_k(x)=\zeta\left(\frac{x-y_k}{R_1}\right),\ \ \rho_k(x)=\rho\left(\frac{x-y_k}{R_k}\right),\
\end{equation*}
where $x\in \mathbb{R},$ and $R_1>R_0$ chosen sufficiently large that
\begin{equation}\label{firstinq}
\left\vert \int_{-\infty}^{\infty} P(\zeta_kf_{j,n_k},\zeta_k g_{n_k})- Q_{1,k}\ dx\right\vert \leq \varepsilon
\end{equation}
and
\begin{equation}\label{2ndInq}
\left\vert \int_{-\infty}^{\infty} P(\rho_kf_{j,n_k},\rho_kg_{n_k})- Q_{2,k}\ dx\right\vert \leq \varepsilon.
\end{equation}
In the last two inequalities we used the notation
\begin{equation*}
P(\phi u_j,\phi v)=|(\phi v)^\prime|^2+|\phi v|^2+\sum_{j=1}^{N}\left(|(\phi u_j)^\prime|^2+|\phi u_j|^2 \right).
\end{equation*}
To see that \eqref{firstinq} and \eqref{2ndInq} are possible, first note that using the first inequality
in \eqref{dichoeq} and the assumptions \eqref{supporteq}, we have that
\begin{equation}\label{techInq}
\begin{aligned}
& \int_{|x-y_k|\leq R_0} \left\vert Q _{n_k}-Q_{1, k}\right\vert\ dx \leq \varepsilon,\\
& \int_{|x-y_k|\geq 2R_k} \left\vert Q _{n_k}-Q_{2, k}\right\vert dx \leq \varepsilon,\
 \int_{A(y_k;R_0, 2R_k)} Q_{n_k}\ dx \leq \varepsilon.
\end{aligned}
\end{equation}
where $A(a;r,R)$ denotes the set $\{x:r\leq |x-a|\leq R \}$ for any $a\in \mathbb{R}, r>0,$ and $R>0.$
The left side of \eqref{firstinq} can be written as
\begin{equation*}
\begin{aligned}
&L= \left\vert \int_{|x-y_k|\leq 2R_1} P(\zeta_kf_{j,n_k},\zeta_k g_{n_k})- Q_{1,k}\ dx\right\vert \\
&= \left\vert \int_{|x-y_k|\leq R_0} Q_{n_k}-Q_{1,k}\ dx \right\vert+\max_{x\in \mathbb{R}}\Omega(R_1;\zeta(x))\int_{A(y_k;R_0, 2R_1)}Q_{n_k}\ dx,
\end{aligned}
\end{equation*}
where for any $a\in \mathbb{R}$ and $\varphi \in C_0^{\infty},$ $\Omega(a;\varphi(x))$ is given by
\begin{equation*}
\Omega(a;\varphi(x))=|\varphi(x)|^2+\frac{1}{a}|\varphi^\prime(x)|^2,\ x\in \mathbb{R}.
\end{equation*}
Using relations \eqref{techInq}, we have that
$L\leq \epsilon + \epsilon = O(\epsilon),$ as $\epsilon\to 0.$ Similarly,
we write the left side of the inequality \eqref{2ndInq} as
\begin{equation*}
\begin{aligned}
& L_1 = \left\vert \int_{|x-y_k|\geq R_k} P(\rho_k f_{j,n_k},\rho_k g_{n_k})- Q_{2,k}\ dx\right\vert \\
&\leq \left\vert \int_{A(y_k;R_k, 2R_k)} P(\rho_k f_{j,n_k},\rho_k g_{n_k})- Q_{2,k}\ dx\right\vert
+ \left\vert \int_{|x-y_k|\geq 2R_k}Q_{n_k}-Q_{2,k}\ dx \right\vert \\
& \leq \max_{x\in \mathbb{R}}\Omega(R_k;\rho(x))\int_{A(y_k; R_k, 2R_k)}Q_{n_k}\ dx
+ \int_{|x-y_k|\geq 2R_k}|Q_{n_k}-Q_{2,k}|\ dx,
\end{aligned}
\end{equation*}
and hence, $L_1 \leq \epsilon + \epsilon = O(\epsilon),$ as $\epsilon\to 0.$
Let us now define $\Delta_{k}^{(1)}$ and $\Delta_{k}^{(2)}$ by setting
\begin{equation*}
\left\{
\begin{aligned}
& \Delta_{k}^{(1)}= (\mathbf{f}_{j,k}^{(1)},g_{k}^{(1)})=(\zeta_k\mathbf{f}_{j,n_k}, \zeta_k g_{n_k})=\zeta_k \Delta_{n_k}, \\
& \Delta_{k}^{(2)}=(\mathbf{f}_{j,k}^{(2)},g_{k}^{(2)})=(\rho_k \mathbf{f}_{j,n_k},\rho_k g_{n_k}) =\rho_k\Delta_{n_k},
\end{aligned}
\right.
\end{equation*}
and let $\Theta_k=(\mathbf{u}_{j,k}, v_{k})$ be such that $\Delta_{n_k}=\Delta_{k}^{(1)}+\Delta_{k}^{(2)}+\Theta_k.$
Then $\Delta_{k}^{(1)}, \Delta_{k}^{(2)}, \Theta_k$ are all in $\mathcal{Y}.$
Since $\int_{-\infty}^{\infty}|F(\Delta_{k}^{(1)})|\ dx $ is bounded, there exists a
subsequence of $\{\Delta_{k}^{(1)}\}_{k\in \mathbb{N}},$ which we denote again by the same symbol,
and a positive real number $\bar{\lambda}=\bar{\lambda}(\epsilon)$ such that
$\int_{-\infty}^{\infty}F(\Delta_{k}^{(1)})\ dx \to \bar{\lambda}.$
Then, for sufficiently large $k,$
\begin{equation}
\int_{-\infty}^{\infty}F(\Delta_{k}^{(1)})\ dx - \bar{\lambda}=O(\epsilon).
\end{equation}
To estimate the proportion of the constraint functional carried by $\Delta_{k}^{(2)},$
we write the integral $\int_{-\infty}^{\infty}F(\Delta_{n_k}) \ dx$ as
\begin{equation}\label{estimate3}
 \int_{-\infty}^{\infty}F(\Delta_{k}^{(1)})+\int_{-\infty}^{\infty}F(\Delta_{k}^{(2)})
 +\int_{A(y_k;R_0, 2R_k)}\left[F(\Delta_{n_k})-F(\Delta_{k}^{(1)})-F(\Delta_{k}^{(1)})\right],
\end{equation}
where all integrals are with respect to $x.$ The last integral in this equation is estimated
as follows:
\begin{equation*}
\begin{aligned}
& \int_{A(y_k;R_0, 2R_k)}\left[F(\Delta_{n_k})-F(\Delta_{k}^{(1)})-F(\Delta_{k}^{(1)})\right]\ dx
\leq C \|\Theta_k\|_{\mathcal{Y}}^{2}\\
&\leq \max\left\{|1-\zeta_k-\eta_k|_{\infty}^{2}, \frac{|\zeta^{\prime}|_{\infty}^{2}}{R_{1}^{2}}+\frac{|\eta^{\prime}|_{\infty}^{2}}{R_{k}^{2}} \right\}
\int_{A(y_k; R_1,2R_k)}Q_{n_k}\ dx=O(\epsilon),
\end{aligned}
\end{equation*}
as $\epsilon\to 0.$ Thus, from \eqref{estimate3}, we can conclude that
\begin{equation*}
\int_{-\infty}^{\infty}F(\Delta_{n_k})\ dx= \int_{-\infty}^{\infty}F(\Delta_{k}^{(1)})\ dx+\int_{-\infty}^{\infty}F(\Delta_{k}^{(2)})\ dx+O(\epsilon).
\end{equation*}
It then follows by taking the limit of both sides as $k\to \infty$ that
\begin{equation*}
\int_{-\infty}^{\infty}F(\Delta_{k}^{(2)})\ dx=\lambda - \bar{\lambda}+O(\epsilon).
\end{equation*}
To prove the assertion that the sum of the values of $\mathsf{K}$
at $\Delta_{k}^{(1)}$ and $\Delta_{k}^{(2)}$ does not exceed $\mathsf{K}(\Delta_{n_k}),$
we write
\begin{equation}\label{Kasse}
\begin{aligned}
\mathsf{K}(\Delta_{n_k})&=\mathsf{K}\left(f_{1,k}^{(1)}+f_{1,k}^{(2)}+u_{1,k},\ .\ .\ . \ , f_{N,k}^{(1)}+f_{N,k}^{(2)}+u_{N,k}, g_k^{(1)}+g_k^{(2)}+v_k\right)\\
& = \mathsf{K}(\mathbf{f}_{j,k}^{(1)}, g_k^{(1)})+\mathsf{K}(\mathbf{f}_{j,k}^{(2)}, g_k^{(2)})+\mathsf{K}(\mathbf{u}_{j,k}, v_k)+\sum_{j=1}^{N}J_j+J,
\end{aligned}
\end{equation}
where the integrals $J$ and $J_j$ on the right-hand side are given by
\begin{equation*}
\begin{aligned}
J =\gamma\int_{-\infty}^{\infty}& \left[(v_k)^\prime (\zeta_kg_{n_k})^\prime +\left(v_k+(\zeta_kg_{n_k})^\prime \right)(\rho_kg_{n_k})^\prime \right]\ dx\\
& +c_\tau \int_{-\infty}^{\infty} \left[v_k\zeta_kg_{n_k}+(v_k+\zeta_kg_{n_k})\rho_kg_{n_k} \right]\ dx
\end{aligned}
\end{equation*}
and for each $j=1,...,N,$
\begin{equation*}
\begin{aligned}
J_j=\int_{-\infty}^{\infty} & \left[(u_{j,k})^\prime (\zeta_kf_{j,n_k})^\prime +\left((u_{j,k})^\prime
+(\zeta_kf_{j,n_k})^\prime \right)(\rho_kf_{j,n_k})^\prime \right]\ dx\\
& + \sigma \int_{-\infty}^{\infty} \left[u_{j,k}\zeta_k f_{j,n_k}+\left(u_{j,k}+\zeta_kf_{j,n_k}\right)\rho_kf_{j,n_k} \right]\ dx.
\end{aligned}
\end{equation*}
From the Cauchy-Schwarz inequality, it follows that
\begin{equation*}
J\leq C\|\Theta_k\|_{\mathcal{Y}}\cdot \|g_{n_k}\|=O(\epsilon)\ \ \textrm{and}\ \ J_j\leq C\|\Theta_k\|_{\mathcal{Y}}\cdot \|f_{j,n_k}\|=O(\epsilon),
\end{equation*}
where $C=C(\zeta_k,\rho_k).$ Thus, from \eqref{Kasse}, we obtain that
\begin{equation*}
\mathsf{K}(\Delta_{n_k})=\mathsf{K}(\Delta_{k}^{(1)})
+\mathsf{K}(\Delta_{k}^{(2)})+O(\epsilon).
\end{equation*}
To complete the proof of Lemma~\ref{Maindicho},
it only remains to establish inequalities in \eqref{Kestimate}. To establish the first inequality, we see that
\begin{equation*}
\begin{aligned}
C & \|\Delta_{k}^{(1)}\|_{\mathcal{Y}}^{2} \geq \mathsf{K}(\Delta_{k}^{(1)})=\mathsf{K}(\zeta_kf_{1,n_k},...,\zeta_kf_{N,n_k},\zeta_kg_{n_k})\\
& = O(\epsilon)+\int_{-\infty}^{\infty}\zeta_k^{2} (f_{1,n_k},\ldots,f_{N,n_k}, g_{n_k})D(L_{ii})(f_{1,n_k},\ldots,f_{N,n_k}, g_{n_k})^T\ dx \\
& \geq O(\epsilon)+\int_{-\infty}^{\infty}\zeta_k^{2} Q_{n_k}\ dx = O(\epsilon)+\int_{B(y_k,R_0)}Q_{n_k}\ dx
+ \int_{A(y_k; R_0, 2R_k)}\zeta_k^{2} Q_{n_k}\ dx\\
& = \int_{-\infty}^{\infty}Q_{1,k}\ dx+ O(\epsilon) \geq \bar{\mu}+O(\epsilon),
\end{aligned}
\end{equation*}
so the first inequality in \eqref{Kestimate} has been proved. The second inequality in \eqref{Kestimate} can be proved similarly. \hfill $\Box$

\medskip

\begin{lem}\label{dicholem2}
Let $Q_n$ be as defined in \eqref{Qdef} and that the dichotomy case occurs.
Then there exists $\theta \in (0,\lambda)$ such that
\begin{equation*}
I_{\lambda} \geq I_{\theta}+ I_{({\lambda-\theta})}.
\end{equation*}
\end{lem}
\noindent {\bf Proof.}
Let $\bar{\lambda}=\bar{\lambda}(\epsilon)$ be as defined in Lemma~\ref{Maindicho}. Since $\int_{-\infty}^{\infty} F(\Delta_{k}^{(1)}) dx$ is bounded,
the range of values of $\bar{\lambda}(\epsilon)$ remains
bounded as $\epsilon \to 0.$ Thus, by restricting attention to a sequence of
values of $\epsilon$ tending to $0$ and extracting an appropriate subsequence from this sequence,
we may assume that $\bar{\lambda}(\epsilon)\to \theta $ as $\epsilon \to 0.$
It is claimed that $\theta \in (0,\lambda).$ To see this, first notice that from
\begin{equation*}
\mathsf{K}(\Delta_{n_k})=\mathsf{K}(\Delta_{k}^{(1)})
+\mathsf{K}(\Delta_{k}^{(2)})+O(\epsilon),
\end{equation*}
it follows immediately that
\begin{equation}\label{limitinf}
I_\lambda= \liminf_{k} \mathsf{K}(\Delta_{n_k})\geq \liminf_{k}\mathsf{K}(\Delta_{k}^{(1)})
+\liminf_{k}\mathsf{K}(\Delta_{k}^{(2)})+O(\epsilon).
\end{equation}
Suppose for the sake of contradiction that $\theta \leq 0.$ Then we have that
\begin{equation*}
\int_{-\infty}^{\infty}F(\Delta_{k}^{(2)})(x)\ dx= \lambda - \theta +O(\epsilon)
\end{equation*}
for sufficiently large $n.$ Let us now define $\bar{\Delta}_{k}^{(2)}=\phi_k \Delta_{k}^{(2)},$ where $\phi_k$ is chosen
such that $\int_{-\infty}^{\infty} F(\bar{\Delta}_{k}^{(2)}) dx=\lambda-\theta.$
Then $\phi_k=1+O(\epsilon)$ and
\begin{equation}\label{EIQ}
\mathsf{K}(\Delta_{k}^{(2)})=\frac{1}{\phi_{k}^{2}}\mathsf{K}(\bar{\Delta}_{k}^{(2)})\geq \frac{1}{\phi_{k}^{2}}I_{\lambda-\theta}
\geq \frac{1}{(1+O(\epsilon))^2}I_\lambda,
\end{equation}
where the last inequality is a consequence of Lemma~\ref{monotono}. From \eqref{limitinf}, \eqref{EIQ}, and
the first inequality of \eqref{Kestimate}, it follows that
\begin{equation*}
I_\lambda \geq C\bar{\mu}
+\frac{1}{(1+O(\epsilon))^2}I_\lambda+O(\epsilon).
\end{equation*}
As $\epsilon\to 0,$ the last inequality gives $I_\lambda \geq C\bar{\mu} + I_\lambda > I_\lambda,$ a contradiction.

\smallskip

\noindent On the other hand, if it were the case that $\theta \geq \lambda,$ then we would have
$\int_{-\infty}^{\infty}F(\Delta_{k}^{(1)})\ dx= \theta +O(\epsilon)$ for sufficiently large $n,$ and a similar argument as
in the case $\theta \leq 0$ would show that \eqref{limitinf} yields
\begin{equation*}
I_\lambda \geq C (\mu-\bar{\mu})
+\frac{1}{(1+O(\epsilon))^2}I_\lambda+O(\epsilon),
\end{equation*}
which implies $I_\lambda \geq C(\mu-\bar{\mu}) + I_\lambda > I_\lambda,$ another contradiction.
This proves the claim that $\theta \in (0,\lambda).$

\smallskip

\noindent Finally, as a consequence the above arguments, one also obtains that
\begin{equation*}
I_\lambda \geq \frac{1}{(1+O(\epsilon))^2}I_\theta
+\frac{1}{(1+O(\epsilon))^2}I_{\lambda-\theta}+O(\epsilon),
\end{equation*}
which upon taking limit as $\epsilon\to 0,$ gives $I_\lambda \geq I_\theta+ I_{\lambda-\theta},$ proving the lemma.\hfill $\Box$

\medskip

We can now rule out the dichotomy condition:

\begin{lem}
The dichotomy does not occur.
\end{lem}
\noindent {\bf Proof.}
This follows from Lemma~\ref{dicholem2} and Corollary~\ref{strictsubadd}. \hfill $\Box$

\medskip

\noindent With both vanishing and dichotomy alternatives ruled out, we can now complete the proof of Theorem~\ref{existence}. Because
vanishing and dichotomy have been ruled out, Lemma~\ref{concentration} guarantees that sequence $%
\{Q_{n}\}$ is tight, i.e., there exists a sequence of real numbers $%
\{y_{n}\}_{n\in \mathbb{N}}$ such that for any $\varepsilon >0,$ there exists $R=R(\varepsilon
)$ so that for all $n\in \mathbb{N},$%
\begin{equation*}
\int_{|x-y_n|\leq R}Q_{n}(x)\ dx \geq \mu-\varepsilon, \ \ \int_{|x-y_n|\geq R}Q_{n}(x)\ dx \leq \varepsilon,
\end{equation*}
and%
\begin{equation*}
\left| \int_{|x-y_n|\geq R}F\left(\textbf{f}_{j,n},g_n\right)\ dx \right|\leq C\ \|(\textbf{f}_{j,n},g_n)\|_{\mathcal{Y}} \int_{|x-y_n|\geq R}Q_{n}(x)\ dx = O(\epsilon),
\end{equation*}
as $\epsilon \to 0.$ It then follows that for $n$ large enough,
\begin{equation}\label{bound1}
\left| \int_{|x-y_n|\leq R} F\left(\textbf{f}_{j,n},g_n\right)\ dx -\lambda \right| \leq \epsilon.
\end{equation}
Denote by $w_{j,n}, 1\leq j \leq N,$ and $z_{n}$ the translated functions
\begin{equation*}
w_{j,n}(x)=f_{j,n}(\cdot +y_{n}),\ \
z_{n}(x)=g_{n}(\cdot +y_{n}).
\end{equation*}
Then, $\{(\mathbf{w}_{j,n},z_n)\}$ is bounded in $\mathcal{Y},$ and hence by the Banach-Alaoglu theorem, there
exists a subsequence, we again label
by $\{(\mathbf{w}_{j,n},z_n)\},$ which converges weakly in $\mathcal{Y}$ to a vector-function $(\Phi_{1},...,\Phi_N, \Psi).$
It then follows immediately from \eqref{bound1} that for $n\geq n_0,$
\begin{equation}\label{bound2}
\lambda \geq \int_{-R }^{R }F\left(w_{1,n},..., w_{N,n}, z_n \right)\ dx \geq \lambda - \epsilon.
\end{equation}
Since $H^1([-R,R])$ is compactly embedded in $L^2([-R,R]),$ we have
\begin{equation*}
\begin{aligned}
\int_{-R }^{R }& \left\vert w_{1,n}^{2}z_{n}- \Phi_{1}^{2} \Psi\right\vert \ dx \leq
|w_{1,n}+\Phi_{1}|_{\infty}\cdot \|z_n\| \cdot \|w_{1,n}-\Phi_{1}\|_{L^2(-R,R)} \\
& \ \ \ \ \ \ \ \ \ \ \ \ \ \ \ \ \ \ \  \ \ \ \ \ \ +\|w_{1,n}\|_1^{2}\cdot \|z_n-\Psi\|_{L^2(-R,R)}\\
&\leq C \left( \|w_{1,n}-\Phi_{1}\|_{L^2(-R,R)}+\|z_n-\Psi \|_{L^2(-R,R)} \right)\to 0,
\end{aligned}
\end{equation*}
as $n\to \infty.$ Similarly, $\int_{-R }^{R }w_{j,n}^{2}z_{n}\ dx \to
\int_{-R }^{R }\Phi_{j}^{2}\Psi \ dx $ for all $ 2 \leq j \leq N.$
We also have
\begin{equation*}
|z_n-\Psi|_{L^3(-R,R)}\leq C \|z_n-\Psi\|_1^{1/6} \|z_n-\Psi\|_{L^2(-R,R)}^{5/6}\leq C \|z_n-\Psi\|_{L^2(-R,R)}^{5/6},
\end{equation*}
and hence, $\int_{-R}^{R}z_n^{3}\ dx \to \int_{-R}^{R}\Psi^{3}\ dx.$
Therefore, from \eqref{bound2}, we have that
\begin{equation*}
\lambda \geq \int_{-R }^{R }F\left(\Phi_{1},...,\Phi_{N},\Psi \right)\ dx \geq \lambda - \epsilon.
\end{equation*}
Thus, for $\epsilon=1/j,\ j\in \mathbb{N},$ there exists $R_j>j$ such that
\begin{equation*}
\lambda \geq \int_{-R_j }^{R_j }F\left(\Phi_{1},...,\Phi_{2},\Psi \right)\ dx \geq \lambda - \frac{1}{j},
\end{equation*}
and consequently, as $j\to \infty,$ we have that $(\Phi_{j}, \Psi)\in \mathcal{A}.$ Furthermore, from the
weak lower semicontinuity of $\mathsf{K}$ and the invariance $\mathsf{K}$ by translations, we have
\begin{equation*}
I_\lambda = \lim_{n\to \infty}\mathsf{K}(\mathbf{f}_{j,n},g_n)\geq \mathsf{K}(\Phi_{j},\Psi)\geq I_\lambda,
\end{equation*}
and thus, $(\Phi_{j},\Psi)$ must be a minimizer for $I_\lambda,$ i.e., $(\Phi_{j},\Psi)\in S(\lambda).$
But then $(\Phi_{j},\Psi)$ must satisfy the Euler-Lagrange equation for (P1), i.e.,
there exists some multiplier $\kappa\in \mathbb{R}$ (Lagrange multiplier) such that
\begin{equation}\label{EulerEq}
\left\{
\begin{aligned}
 &-\Phi_{1}^{\prime \prime}+\sigma \Phi_{1}=\kappa \alpha_1 \Psi \Phi_{1}, \\
 &\ \ \ \ \ldots \ \ \ \  \ \ \ \ \ldots \ \ \ \  \ \  \ldots \\
 &-\Phi_{N}^{\prime \prime}+\sigma \Phi_{N}=\kappa \alpha_N \Psi \Phi_{N},  \\
& -\gamma \Psi^{\prime \prime}+ c_\tau \Psi =\frac{\kappa}{2} \left(\beta \Psi^2 +\alpha_1 \Phi_{1}^{2}+\ldots +\alpha_1 \Phi_{N}^{2} \right).
\end{aligned}
\right.
\end{equation}
An easy calculation proves that the Lagrange multiplier is positive:
\begin{prop}
The Lagrange multiplier satisfies $\kappa>0.$
\end{prop}
\noindent{\bf Proof.}
Multiplying the first and second equations above by
$\Phi_{j}$ and $\Psi,$ respectively, and integrating over the real line, we obtain
\begin{equation*}
\left\{
\begin{aligned}
& \int_{-\infty}^{\infty}\left( (\Phi_{j}^\prime)^2 + \sigma \Phi_{j}^{2} \right) \ dx
= k \int_{-\infty}^{\infty} \alpha_j\ \Psi \Phi_{j}^2\ dx,\ j=1,2,...,N, \\
& \int_{-\infty}^{\infty}\left( \gamma (\Psi^\prime)^2+c_\tau \Psi^2 \right) \ dx = \frac{\kappa}{2}\int_{-\infty}^{\infty} \left(\beta\ \Psi^3+\sum_{j=1}^{N} \alpha_j\ \Phi_{j}^2 \Psi \right)\ dx.
\end{aligned}
\right.
\end{equation*}
Adding these $N+1$ equations and using the facts that $\mathsf{K}(\Phi_{j},\Psi)=I_\lambda$
and $\int_{-\infty}^{\infty}F(\Phi_{j},\Psi)\ dx = \lambda,$ we obtain
\begin{equation*}
\kappa = \frac{2}{3\lambda} I_\lambda > 0,
\end{equation*}
which is the desired result. \hfill $\Box$

\medskip

\noindent Finally, we see that
these equations \eqref{EulerEq} are satisfied by $\Phi_{j}, 1\leq j \leq N,$ and $\Psi$ if and only if the functions $u_j$ and $v$ defined by
\begin{equation*}
u_j(x,t) =\kappa e^{i\omega t}e^{ic(x-ct)/2}\Phi_{j}(x-ct),\  v(x,t)= \kappa \Psi(x-ct)
\end{equation*}
are solutions of \eqref{lwsw}. That is, solutions to the variational problem (P1) corresponds to solitary-wave profiles of \eqref{lwsw}.

\smallskip

\noindent To complete the proof of Theorem~\ref{existence}, it
only remains to prove smoothness and exponential decay of the solutions:
\begin{prop}\label{smoothdecay}
Suppose $(\Phi_1,...,\Phi_N,\Psi)\in \mathcal{Y} $ is a solution of \eqref{ODE}, in the sense of distributions. Then
\begin{itemize}
\item[(i)] $\Phi_1,...,\Phi_N,\Psi \in H^{\infty}(\mathbb{R}).$
\item[(ii)] One has pointwise exponential decay, i.e.,
\begin{equation*}
|\Phi_j(x)|\leq Ce^{-\delta_j |x|},\ \textrm{and} \ |\Psi(x)|\leq Ce^{-\delta |x|},
\end{equation*}
holds for all $x\in \mathbb{R},$ where $\delta, \delta_j>0,$ and $C>0$ are suitable constants.
\end{itemize}
\end{prop}
\noindent {\bf Proof.}
Statement (i) follows by a standard bootstrap argument.
Since $\Phi_1, ..., \Phi_N,$ $\Psi$ are $H^1$ functions,
and $H^1$ is an algebra, it follows that $\Psi^2, \Phi_j^{2},$ and $\Phi_j\Psi, 1\leq j \leq N,$ are also $H^1$ functions. Since the convolution
operation with $K_s$ takes $H^s$ to $H^{s+2}$ for any $s\geq 0,$ so \eqref{convrepre} implies that $\Phi_j,...,\Phi_N, \Psi$ are in $H^3.$
But then $\Psi^2, \Phi_j^{2},$ and $\Phi_j\Psi$ are $H^3$ functions, so \eqref{convrepre} implies
that $\Phi_j$ and $\Psi$ are in $H^5,$ and so on. Continuing this argument inductively proves that $\Phi_1, ..., \Phi_N,$ $\Psi$ are in $H^\infty.$

\smallskip

\noindent To prove decay estimates, we borrow an argument from the proof of Theorem~8.1.1 of \cite{[Ca]}.
Fix $j\in \{1,2,...,N\}.$ For $\epsilon>0$ and $\delta>0,$ consider the
function $ \varphi(x)=e^{\epsilon |x|/(1+\delta |x|)}\in L^\infty(\mathbb{R}).$
Multiplying the first equation in \eqref{ODE} by $\varphi \Phi_j,$ we get
\begin{equation*}
-\int_{-\infty}^{\infty} \varphi\ \Phi_j^{\prime \prime} \Phi_j\ dx + \sigma \int_{-\infty}^{\infty}\varphi\ \Phi_j^2\ dx = \alpha_j\int_{-\infty}^{\infty}\varphi\ \Phi_j^2 \Psi \ dx.
\end{equation*}
Integrating by parts and using the fact that $\varphi^\prime\leq \epsilon \varphi,$ we get
\begin{equation*}
\sigma \int_{-\infty}^{\infty}\varphi \Phi_j^2\ dx\leq \int_{-\infty}^{\infty} \varphi (\Phi_j^\prime)^2\ dx + \epsilon \int_{-\infty}^{\infty} \varphi |\Phi_j\Phi_j^\prime|\ dx
+ \alpha_j\int_{-\infty}^{\infty}\varphi \Phi_j^2 | \Psi| \ dx.
\end{equation*}
Now using the Cauchy-Schwarz inequality, we obtain from the preceding inequality that
\begin{equation*}
\begin{aligned}
\left(\sigma-\frac{\epsilon}{2} \right)\int_{-\infty}^{\infty}\varphi \Phi_j^2\ dx & \leq \left(1+ \frac{\epsilon}{2} \right)\int_{-\infty}^{\infty}\varphi (\Phi_j^\prime)^2\ dx + \alpha_j\int_{-\infty}^{\infty}\varphi \Phi_j^2 | \Psi| \ dx \\
& \leq \alpha_j\int_{-\infty}^{\infty}\varphi \Phi_j^2 | \Psi| \ dx,
\end{aligned}
\end{equation*}
with $\epsilon$ chosen to be sufficiently small. Thus, for $\epsilon$ small enough, we deduce that
\begin{equation}\label{decay1}
\int_{-\infty}^{\infty}\varphi(x) \Phi_j^2(x)\ dx \leq C \int_{-\infty}^{\infty}\varphi(x) \Phi_j^2(x) | \Psi(x)| \ dx,
\end{equation}
where $C=C(\epsilon, \sigma, \alpha_j)$ (independent of $\delta$).
Since $\Psi$ is an $H^1$ function, then $\Psi(x)\to 0$ as $|x|\to \infty.$ We can find $R>0$ such that $|\Psi(x)|\leq 1/(2C)$ for $|x|\geq R.$
It then follows from \eqref{decay1} that
\begin{equation*}
\int_{-\infty}^{\infty}\varphi(x) \Phi_j^2(x)\ dx \leq 2C \int_{B_R} \varphi(x) \Phi_j^2(x) | \Psi(x)| \ dx.
\end{equation*}
Taking $\delta\to 0,$ Fatou's lemma and Lebesgue's theorem yields
\begin{equation}\label{decay2}
\int_{-\infty}^{\infty} e^{\epsilon |x|} |\Phi_j(x)|^2 \ dx< \infty,\ j=1,2,...,N.
\end{equation}
Now since $\Phi_j$ belongs to $H^1,$ then $\Phi_j(x)\to 0$ as $|x|\to \infty$ and $\Phi_j$ is
globally Lipschitz continuous on $\mathbb{R}.$ From these two
properties of $\Phi_j$ and \eqref{decay2}, one can easily show that $e^{\delta_1 |x|}\Phi_j(x)$ is bounded on $\mathbb{R}$ for some $0<\delta_1 \leq \epsilon$ (for details, the reader may consult the proof of Theorem~8.1.1 in \cite{[Ca]}).

\smallskip

\noindent To obtain the decay estimate for $\Psi,$ multiplying the second equation in \eqref{ODE} by $\varphi \Psi$ gives, as above, the following estimate
\begin{equation*}
\int_{-\infty}^{\infty}\varphi\ \Psi^2\ dx \leq C \int_{-\infty}^{\infty} \varphi\ \left(|\Psi|^3+\sum_{j=1}^{N} \Phi_j^{2}|\Psi| \right)\ dx,
\end{equation*}
where $C=C(\epsilon, \gamma, c, \tau, \beta, \alpha_j).$ Choose $\epsilon < 2\delta_1.$
Then, using decay estimates for the functions $\Phi_1,...,\Phi_N$ proved above, we can show,
as before, that $\int_{-\infty}^{\infty} \varphi \Psi^2\ dx$ is
bounded by some constant which is independent of $\delta.$ Then, taking $\delta\to 0.$ allows us to deduce that
\begin{equation*}
\int_{-\infty}^{\infty} e^{\epsilon |x|} |\Psi(x)|^2 \ dx< \infty,
\end{equation*}
and from here we can proceed as we did for $\Phi_j(x).$\hfill $\Box$

\smallskip

The proof of Theorem~\ref{existence} is now complete.

\section{Properties of Solitary Waves}\label{Section3}

\noindent In this section we establish some properties of travelling solitary waves.
To do so, we take advantage of the convolution representation for
solutions of the equation \eqref{ODE}. We assume throughout this section, unless otherwise stated,
that the assumptions \eqref{assumptions} hold with $\gamma >0.$

\noindent Provided  $\eta=c_\tau/\gamma>0$ in \eqref{ODE},
we can rewrite \eqref{ODE} in the form
\begin{equation}\label{convrepre}
\left\{
\begin{aligned}
& \Phi_j= \alpha_j\ K_{\sigma}\star \Phi_j\Psi,\ \ 1\leq j \leq N, \\
& \Psi= \frac{1}{2\gamma} \ K_{\eta}\star \left( \beta \Psi^{2}+ \sum_{j=1}^{N}\alpha_j \Phi_j^{2} \right),
\end{aligned}
\right.
\end{equation}
where for any $s>0,$ the kernel $K_{s}$ is defined explicitly in terms of its Fourier symbol
\begin{equation}\label{qdef}
q(\xi)=\widehat{K}_s(\xi)=\frac{1}{s+\xi^2}.
\end{equation}
The kernel $K_{s}$ defined via its Fourier transform as in \eqref{qdef} is a real-valued, even, bounded,
continuous function, and $K_{s}(x)\to 0$ as $|x|\to \infty.$
Furthermore, $K_{s}$ is strictly positive on $\mathbb{R}.$
To see this, one can use the Residue Theorem and Jordan's
lemma (see, for example, Chapter 3 of \cite{[Stein]}) to represent $K_{s}$ explicitly in the form
\begin{equation}\label{KsDef}
K_s(x)=\frac{1}{2\sqrt{s}}\sqrt{\frac{\pi}{2}}e^{-\sqrt{s}|x|}, \ s>0, \ x\in \mathbb{R}.
\end{equation}

The first property concerns signs of $\Phi_j, 1\leq j \leq N,$ and $\Psi:$

\begin{thm}\label{signthm}
Every solution $(\Phi_1,...,\Phi_N,\Psi)$ of the \eqref{ODE}
satisfies the following properties:
\begin{itemize}
\item[(i)] $\Psi(x)> 0$ for all $x\in \mathbb{R}.$
\item[(ii)] The functions $\Phi_1(x),...,\Phi_N(x)$ are of one sign on $\mathbb{R}.$
\end{itemize}
\end{thm}

\noindent {\bf Proof.} To prove (i), we use the convolution representation of $\Psi,$ namely
\begin{equation}\label{Psiconvo}
\Psi(x)=\frac{1}{2\gamma} \int_{-\infty}^{\infty}K_{\eta}(x-y)Q(y)\ dy \ \textrm{with}\ Q(y)= \beta \Psi^{2}(y)+ \sum_{j=1}^{N}\alpha_j\ \Phi_j^{2}(y).
\end{equation}
Since the kernel $K_\eta$ is a strictly positive and $Q(x)$ is everywhere non-negative, it then
follows from \eqref{Psiconvo} that $\Psi(x)>0$ everywhere provided that $(\Phi_1,...,\Phi_N,\Psi)$
is a solution of \eqref{ODE} with at least one component being
nonzero on a set of positive measure.

\smallskip

\noindent To prove (ii), denote $\Theta=(\Phi_1,...,\Phi_N, \Psi).$ Let $U_j=|\Phi_j|, j=1,..,N,$ and denote $\Delta=(U_1,...,U_N,\Psi).$
It is a standard fact from real analysis that if $\Phi_j\in H^1,$ then $|\Phi_j(x)|$ is in $H^1$ and
\begin{equation}\label{Dnorm}
\int_{-\infty}^{\infty}||\Phi_j|_x|^2\ dx \leq \int_{-\infty}^{\infty} |(\Phi_j)_x|^2\ dx.
\end{equation}
(For a proof of this elementary fact, readers may consult Lemma~3.4 of \cite{[AB2]}.)
Then $\Delta\in \mathcal{Y}$ and using \eqref{Dnorm}, it follows that $\mathsf{K}(\Delta)\leq \mathsf{K}(\Theta).$
Since $\int_{-\infty}^{\infty}F(\Delta)\ dx = \int_{-\infty}^{\infty}F(\Theta)\ dx,$ thus $\Delta$ and $\Theta$ are both
in $S(\lambda)$ for some $\lambda>0.$
Observe that since $\Theta$ and $\Delta$ satisfy the same
equations \eqref{EulerEq}, we have that for each $j=1,...,N,$
\begin{equation}\label{Onesign}
\left\{
\begin{aligned}
& -\Phi_{j}^{\prime \prime}(x)+\sigma \Phi_j(x)=\kappa \alpha_j\ \Psi(x) \Phi_j(x), \\
& -U_{j}^{\prime \prime}(x)+\sigma U_j(x)=\kappa \alpha_j\ \Psi(x) U_j(x). \\
\end{aligned}
\right.
\end{equation}
Multiplying the first equation in \eqref{Onesign} by $U_j$ and the second by $\Phi_j$ and
subtracting the second from the first, we see that the Wronskian
\begin{equation*}
 W(\Phi_j(x),U_j(x))=\left\vert
\protect\begin{array}{cc}
\Phi _{j}(x) & U_{j}(x) \protect \\
\Phi _{j}^{\prime }(x) & U_{j}^{\prime }(x) %
\protect\end{array}%
\right\vert = \textrm{constant}.
\end{equation*}
But since $W(\Phi_j,U_j)\to 0$ as $x\to \infty,$ we must have $W(\Phi_j,U_j)=0$ for all $x\in \mathbb{R}.$
Then the functions $\Phi_j$ and $U_j, 1\leq j \leq N,$ are linearly dependent and so, $\Phi_j, 1\leq j \leq N,$ must be of one sign on $\mathbb{R}.$
This completes the proof. \hfill $\Box$

\medskip

\noindent We now use the theory of symmetric decreasing rearrangement
to prove the existence of a solution $(\Phi_{1},..., \Phi_{N},\Psi)$ of \eqref{ODE} such that
$\Phi_{1},...,\Phi_{N}, \Psi$ are even, strictly positive, and decreasing functions on $(0,\infty).$
Recall that, for a non-negative function $w:\mathbb{R}\to [0,\infty),$ one may define the symmetric rearrangement
of $w$ to be the unique function $w^\ast$ with domain $\mathbb{R}$ which has the same
distribution function as $w,$ that is, for every $a>0,$ the sets $\{x:|w(x)|>a\}$ and $\{x:|w^\ast(x)|>a\}$ have the same measure.
In formulas,
\begin{equation*}
w^\ast(x)=\inf \left\{a>0:\frac{1}{2}m(w,a)\leq |x|\right\}=\sup\left\{a>0:\frac{1}{2}m(w,a)> |x|\right\},
\end{equation*}
where $m(w,a)$ denotes the measure of $\{x:|w(x)|>a\}$ for all $a>0$
(or see Chapter~2 of \cite{[Kaw]} for a slightly different but equivalent
definition and also a comprehensive discussion of many different types of rearrangements).
The function $w^\ast$ is clearly radially symmetric and non-increasing in the variable $|x|,$ i.e., $w^\ast(x)=w^\ast(y)$ if $|x|=|y|$ and
$w^\ast(x)\geq w^\ast(y)$ if $|x|\leq |y|.$

\smallskip

\noindent The next two theorems about symmetric decreasing rearrangements play a crucial role in the rest of the paper.

\begin{thm} \label{symprop}
The following statements hold:
\begin{itemize}
\item[(i)] Rearrangement preserves $L^p$ norms, i.e., for every nonnegative function $f$ in $L^p,$
\begin{equation*}
|f|_p=|f^\ast|_p,\ 1\leq p \leq \infty.
\end{equation*}
\item[(ii)] (Hardy-Littlewood Inequality) If $f$ and $g$ are nonnegative measurable functions that vanish at infinity, then
\begin{equation*}
\int_{-\infty}^{\infty}f(x)g(x)\ dx \leq \int_{-\infty}^{\infty}f^\ast (x)g^\ast(x)\ dx.
\end{equation*}
\item[(iii)] (P\'{o}lya-Szeg\H{o} Inequality) The symmetric decreasing rearrangement diminishes $L^2$ norm of the gradient of a positive function $f$ in $H^1:$
    \begin{equation*}
    \int_{-\infty}^{\infty} |(u^\ast)_x|^2\ dx \leq  \int_{-\infty}^{\infty} |u_x|^2\ dx.
    \end{equation*}
\end{itemize}
\end{thm}
\noindent For proofs of these statements, as well as other basic
facts about rearrangements, reader may consult, for example,
the appendix of \cite{[AP]}.

\begin{thm}[F.~Riesz]\label{Riesz}
Let $f_1,....,f_N$ be measurable functions on $\mathbb{R}$ such that $ m\{x:f_j\geq y \} < \infty $ for all $y>0$ and all $1 \leq j \leq N.$ Then
\begin{equation}
|(f_1\star f_2 \star ... \star f_N)(0)|\leq \left[(f_1^\ast)\star (f_2^\ast)\star . . . \star (f_N^\ast)\right](0)
\end{equation}
in the sense that if the right-hand side is finite, then the left-hand side exists and the inequality holds.
\end{thm}
\noindent A proof of Theorem~\ref{Riesz} for $N=3,$ along with a sketch of the inductive proof for $N\geq 3,$ can be found in \cite{[Olver]}.

\smallskip

\noindent Following Weinstein \cite{[Wein]}, we introduce a functional $\Lambda(\Theta)$ by
\begin{equation*}
\Lambda(\Theta)=\frac{\mathsf{K}(\Theta)}{\left(\int_{-\infty}^{\infty}F(\Theta)\ dx \right)^{2/3}},\ \ \Theta=(f_1,...,f_N,g),
\end{equation*}
where $f_1,...,f_N,g\in H^1.$ If $\Lambda(\Theta)$ has a critical point at $\Theta=(\phi_1,...,\phi_N,\psi),$ then
a computation of the Fr\'{e}chet derivative
of $\Lambda$ at $(\phi_1,...,\phi_N,\psi)$ shows
that $(\phi_1,...,\phi_N,\psi)$ is, up to a constant multiple, a solution of \eqref{ODE}.
Consider now the
following unconstrained minimization problem
\begin{equation*}\tag{P2}
\min \left\{ \Lambda(f_1,...,f_N,g): \vec{0} \neq (f_1,...,f_N,g) \in \mathcal{Y} \right\}.
\end{equation*}
\begin{prop}
The problem (P1) is equivalent to (P2).
More precisely, any solution of (P1) is a minimizer of (P2), and if $\Delta=(\phi_1,...,\phi_N,\psi)$ is a minimizer of (P2) then the rescaling
\begin{equation*}
(\phi_1,...,\phi_N,\psi) \mapsto \frac{\lambda^{1/3}}{\left(\int_{-\infty}^{\infty}F(\Delta)\ dx \right)^{1/3}}\ (\phi_1,...,\phi_N,\psi)
\end{equation*}
is a solution of the problem (P1). Moreover, if we define $\Theta$ by
\begin{equation*}
\Theta = \frac{2}{3} \Lambda(\Delta)\ \frac{\Delta}{\left(\int_{-\infty}^{\infty}F(\Delta)\ dx \right)^{1/3}} = \frac{2}{3} I_1\ \frac{U}{ \lambda^{1/3}},
\end{equation*}
then $\Theta$ does not depend on $\lambda$ and solves \eqref{ODE}, and therefore when
substituted into \eqref{SO} yields a travelling solitary wave for \eqref{lwsw}.
\end{prop}
\noindent{\bf Proof.}
The proof follows from Theorem~\ref{existence}.\hfill $\Box$

\medskip

\begin{thm}\label{evenSOL}
Suppose $c_\tau >0$ and $\sigma>0.$
There exists a non-trivial,
solitary-wave solution $(\Phi_1,...,\Phi_N,\Psi)$ of the system \eqref{ODE} such that
\begin{equation}\label{NONneg}
\Phi_1(x)\geq 0,\ . \ . \ ., \Phi_N(x)\geq 0, \Psi(x)\geq 0,
\end{equation}
for $x \in \mathbb{R}.$ Moreover, the functions
$\Phi_1(x),...,\Phi_N(x), \Psi(x)$ can be chosen to be
even, strictly positive, and non-increasing for $x\geq0.$
\end{thm}
\noindent {\bf Proof.}
Let $\Theta=(f_{01},...,f_{0N},g_0)$ be a minimizer of (P2) in $\mathcal{Y}.$
Define $\phi_j=|f_{0j}|, j=1,...,N,$ and $\psi=|g_0|.$ Then $\phi_j$ and $\psi$ are in $H^1,$
and a similar inequality of the form \eqref{Dnorm} holds for $\phi_j$ and $\psi.$
Denote $\Delta=(\phi_1,...,\phi_N,\psi).$ Then $\mathsf{K}(\Delta)\leq \mathsf{K}(\Theta)$ and
\begin{equation*}
\int_{-\infty}^{\infty} F(\Delta)\ dx \geq \int_{-\infty}^{\infty} F(\Theta)\ dx.
\end{equation*}
It follows that $\Lambda(\Delta)\leq \Lambda(\Theta).$ Hence
$\Delta$ is also a minimizer of (P2). Consequently, there exists
a constant $a>0$ such that
such that $(\Phi_1,...,\Phi_N,\Psi)$ defined by
\begin{equation*}
\Phi_j=a \phi_1=a |f_{0j}|,\ 1\leq j \leq N,\ \textrm{and}\ \Psi = a \psi = a |g_0|,
\end{equation*}
is a solution of \eqref{ODE} and \eqref{NONneg} holds.
To prove $\Phi_j>0,$ observe that \eqref{EulerEq} implies $\alpha_j K_\sigma \star (\Psi\Phi_j)=\Phi_j$ for all $1\leq j\leq N,$
where $K_\sigma$ is defined as in \eqref{KsDef}.
It follows that $\Phi_j>0$ on $\mathbb{R}$ for each $1\leq j\leq N.$

\smallskip

\noindent Next, as above, let $\Theta=(f_{01},...,f_{0N},g_0)$ be a minimizer of the problem (P2).
Define $\Delta_\ast=(\phi_{\ast 1},...,\phi_{\ast N},\psi_\ast)$ by
setting $\phi_{\ast j}=f_{0j}^\ast, 1\leq j \leq N,$ and $\psi_\ast=g_0^\ast.$
Then, using the Hardy-Littlewood inequality (Theorem~\ref{symprop}(ii)), we have that
\begin{equation*}
\int_{-\infty}^{\infty}F(\Theta)\ dx\leq \int_{-\infty}^{\infty}F(\Delta_\ast)\ dx.
\end{equation*}
Using the P\'{o}lya-Szeg\H{o} inequality (Theorem~\ref{symprop}(iii)) and
the fact that rearrangement preserves $L^p$ norm (Theorem~\ref{symprop}(i)), it follows
that $\mathsf{K}(\Delta_\ast)\leq \mathsf{K}(\Theta).$ Hence $\Delta_\ast$ is also a minimizer of (P2). Then, as in
the preceding paragraph, there exists $a \in \mathbb{R}$ such that $(\Phi_j,\Psi)=a \Delta_\ast$ is a solution of \eqref{ODE}.
Since $\Phi_j=a \phi_{\ast j}, 1\leq j \leq N,$ and $\Psi=a \psi_\ast$ are
non-increasing functions of $|x|,$ this completes the proof. \hfill $\Box$

\medskip

\noindent We now prove positivity of Fourier transforms of the solitary waves.
\begin{thm}\label{positiveFT}
Suppose $c_\tau >0$ and $\sigma>0.$ Then there exists a solution $(\Phi_1,...,\Phi_N,\Psi)$ of \eqref{ODE} such that
\begin{equation}\label{poFTInq}
\widehat{\Phi}_1(\xi)\geq 0,\ . \ . \ .,\widehat{\Phi}_N(\xi)\geq 0,  \widehat{\Psi}(\xi)\geq 0,
\end{equation}
for $\xi\in \mathbb{R}.$ Moreover,
$\Phi_1,...,\Phi_N, \Psi$ may be chosen so that $\widehat{\Phi}_1,...,\widehat{\Phi}_N, \widehat{\Psi}$
are even, strictly positive, and non-increasing functions of $|\xi|.$
It also follows that $\Phi_1,...,\Phi_N, \Psi$ are even functions.
\end{thm}
\noindent {\bf Proof.}
We follow ideas of Albert \cite{[Albp]}.
Let $\Theta=(f_{01},...,f_{0N},g_0)$ be a minimizer of the problem (P2) in $\mathcal{Y}.$ Choose $u_1,...,u_N,v\in L^2$ such that
\begin{equation*}
\widehat{u}_j=|\widehat{f}_{0j}|\ \ \textrm{and}\ \widehat{v}=|\widehat{g}_0|,\ j=1,2,...,N.
\end{equation*}
The functions $u_1,...,u_N,v$ are real-valued since $|\widehat{f}_{0j}|$
and $|\widehat{g}_0|$ are real-valued and even; and $U=(u_1,...,u_n,v)$ belongs to $\mathcal{Y}.$
Then $\mathsf{K}(U)=\mathsf{K}(\Theta)$ and using Theorem~\ref{Riesz}, we have
\begin{equation*}
\begin{aligned}
\int_{-\infty}^{\infty}F(U)\ dx & = \frac{1}{3}\ \beta\ \widehat{v^3}(0)+\sum_{j=1}^{N}\alpha_j\ \widehat{u_{j}^{2}v}(0)\\
& \geq \frac{1}{3}\ \beta\ \widehat{g_{0}^3}(0)+\sum_{j=1}^{N}\alpha_j\ \widehat{f_{0j}^{2}g_0}(0)=\int_{-\infty}^{\infty}F(\Theta)\ dx.
\end{aligned}
\end{equation*}
It follows that
$\Lambda(U)\leq \Lambda(\Theta).$
Hence $U$ is also a minimizer of (P2). As noted before, there exists a constant $a>0$ such that
\begin{equation*}
(\Phi_1,...,\Phi_N,\Psi)=a(u_1,...,u_N,v)
\end{equation*}
is a solution of \eqref{ODE}. Since $\widehat{\Phi}_j=a|\widehat{f}_{0j}|\geq 0 $
and $\widehat{\Psi}=a|\widehat{g}_0|\geq 0,$ the first assertion of the Theorem follows.

\smallskip

\noindent On the other hand, since for each $1\leq j \leq N,$ $f_{0j}$ is real-valued,
\begin{equation*}
\begin{aligned}
\widehat{\Phi}_j(\xi)& =a|\widehat{f}_{0j}(\xi)| =
\frac{a}{\sqrt{2 \pi}}\left\vert \int_{-\infty}^{\infty}f_{0j}(x)\cos(x\xi)\ dx - i\int_{-\infty}^{\infty}f_{0j}(x)\sin(x\xi)\ dx \right\vert \\
& = \frac{a}{ \sqrt{2\pi}} \sqrt{\left\vert \int_{-\infty}^{\infty}f_{0j}(x)\cos(x\xi)\ dx \right\vert^2+\left\vert\int_{-\infty}^{\infty}f_{0j}(x)\sin(x\xi)\ dx \right\vert^2 }
\end{aligned}
\end{equation*}
is an even function, and hence for all $1\leq j \leq N,$
\begin{equation*}
\begin{aligned}
\Phi_j(x) & =\frac{1}{\sqrt{2\pi}}\int_{-\infty}^{\infty}\widehat{\Phi}_j(\xi)\cos(x\xi) \ d\xi
+\frac{i}{\sqrt{2\pi}}\int_{-\infty}^{\infty}\widehat{\Phi}_j(\xi)\sin(x\xi) \ d\xi \\
& = \frac{1}{\sqrt{2\pi}}\int_{-\infty}^{\infty}\widehat{\Phi}_j(\xi) \cos(x\xi) \ d\xi
\end{aligned}
\end{equation*}
is an even function. Similarly, one can show that $\Psi$ is even.

\smallskip

\noindent To prove positivity of the Fourier transforms of the solitary waves, let
$\Theta=(f_{01},...,f_{0N},g_0)$ be, as above,
a minimizer of (P2). Choose $f_1,...,f_N,g\in L^2$ such that
$\widehat{f}_j=\widehat{f}_{0j}^\ast,\ \textrm{and}\ \widehat{g}=\widehat{g}_{0}^\ast.$ Denote $\Delta=(f_1,...,f_N,g).$
We claim that $\Lambda(\Delta)\leq \Lambda(\Theta).$ First,
using statements (i) and (iii) of Theorem~\ref{symprop}, it follows that $\mathsf{K}(\Delta)\leq \mathsf{K}(\Theta).$
Hence to prove the claim it suffices to show that the denominator of $\Lambda(\Delta)$ is greater than or equal to that of $\Lambda(\Theta).$
Using the inequality of F.~Riesz concerning convolutions of symmetric rearrangements of functions (Theorem~\ref{Riesz}), we have
\begin{equation*}
\begin{aligned}
\int_{-\infty}^{\infty}& F(\Delta) \ dx =\frac{\beta}{3}\left[\widehat{g}\star \widehat{g} \star \widehat{g} \right](0)+\sum_{j=1}^{N}\alpha_j\left[\widehat{f}_j\star \widehat{f}_j \star \widehat{g}\right](0) \\
& = \frac{\beta}{3}\left[\widehat{g}_0^\ast \star \widehat{g}_0^\ast \star \widehat{g}_0^\ast \right](0)+\sum_{j=1}^{N}\alpha_j \left[\widehat{f}_{0j}^\ast \star \widehat{f}_{0j}^\ast \star \widehat{g}_0^\ast \right](0) \\
& \geq \frac{\beta}{3}\left[\widehat{g}_0 \star \widehat{g}_0 \star \widehat{g}_0 \right](0)+\sum_{j=1}^{N}\alpha_j \left[\widehat{f}_{0j} \star \widehat{f}_{0j} \star \widehat{g}_0 \right](0)=\int_{-\infty}^{\infty}F(\Theta)\ dx.
\end{aligned}
\end{equation*}
This proves the claim.
Hence $\Delta$ is also a minimizer of (P2), and it follows that there exists $a\in \mathbb{R}$ such that
\begin{equation*}
(\Phi_{1},...,\Phi_{N},\Psi)=a(f_1,...,f_N,g)
\end{equation*}
is a solution of \eqref{ODE}.
Since $\widehat{\Phi}_j= a \widehat{f}_{0j}^\ast$ and $\widehat{\Psi}=a \widehat{g}_0^\ast$ are
non-increasing functions of $|\xi|,$ it remains
only to show that $\widehat{\Phi}_1,....,\widehat{\Phi}_N,$ and $\widehat{\Psi}$ are everywhere positive.
We only prove that $\widehat{\Phi}_1$ is everywhere positive.
If this is not the case, then
the support of $\widehat{\Phi}_1$ is a finite closed interval $[-a_1,a_1].$ On the other hand,
the support of $\widehat{\Phi}_1\star \widehat{\Psi}$ strictly contains $[-a_1,a_1],$
so that $(\xi^2+\sigma)\widehat{\Phi}_1$ can not equal to $\widehat{\Psi}\star \alpha_1\widehat{\Phi}_1.$
This then contradicts the first equation in \eqref{convrepre} with $j=1.$
Similarly, one can prove that $\Phi_2,...,\Phi_N, \widehat{\Psi}$ are everywhere positive. \hfill $\Box$

\bigskip


\end{document}